\documentclass[11pt]{article}

 \usepackage{setspace}
\usepackage{extdash}

\usepackage{latexsym}
\usepackage{calc}

\usepackage{amssymb}
\usepackage{amsmath}
\usepackage{graphicx}
\usepackage{comment}
\usepackage{psfrag}

\usepackage{float}

\usepackage{color}

\usepackage{nicefrac}

\usepackage{ntheorem}

\usepackage{hyperref}

\usepackage{tikz}
\usetikzlibrary{backgrounds}

\newcounter{hours}\newcounter{minutes}

\def\nr{\par }

\def\beq{\begin{equation}}
\def\eeq{\end{equation}}

\newtheorem{theorem}{Theorem}

\newtheorem{lemma}{Lemma}
\newtheorem{corollary}{Corollary}

\newtheorem{proposition}{Proposition}
\newtheorem{assumption}{Assumption}
\newtheorem{definition}{Definition}

\newtheorem{example}{Example}
\newtheorem{remark}{Remark}
\newcommand{\proof}{\bf Proof: \rm \nr}

\newcommand{\qed}{\hfill $\Box$ \nr \medskip}

\def\ba{\begin{array}}
\def\ea{\end{array}}
\def\beann{\begin{eqnarray*}}
\def\eeann{\end{eqnarray*}}
\def\bea{\begin{eqnarray}}
\def\eea{\end{eqnarray}}

\def\BT{\begin{theorem}}
\def\ET{\end{theorem}}
\def\BL{\begin{lemma}}
\def\EL{\end{lemma}}
\def\BC{\begin{corollary}}
\def\EC{\end{corollary}}
\def\BE{\begin{example}}
\def\EE{\end{example}}
\def\BD{\begin{definition}}
\def\ED{\end{definition}}
\def\BR{\begin{remark}}
\def\ER{\end{remark}}
\def\BAS{\begin{assumption}}
\def\EAS{\end{assumption}}
\def\BI{\begin{itemize}}
\def\EI{\end{itemize}}

\def\BMP{\begin{minipage}{9.5cm}}
\def\EMP{\end{minipage}}
\def\MPT{\begin{minipage}{11.5cm}}
\def\EPT{\end{minipage}}

\def\R{\mathbb{R}}

\newcommand*\samethanks[1][\value{footnote}]{\footnotemark[#1]}
\newcommand{\NNorm}[2]{\left\Vert {#1} \right\Vert_{#2}}

\title{Global aspects of the continuous reformulation  \\ for cardinality-constrained optimization problems}
\author{
S. L\"ammel
\thanks{
Department of Mathematics, Chemnitz University of Technology,
Reichenhainer Str. 41, 09126
Chemnitz, Germany; e-mail: sebastian.laemmel@mathematik.tu-chemnitz.de, vladimir.shikhman@mathematik.tu-chemnitz.de.
 } \and V. Shikhman\samethanks[1]
}

\begin{document}
\thispagestyle{empty}
\maketitle
\vspace{-5ex}
\abstract{
 The main goal of this paper is to relate the topologically relevant stationary points of a cardinality-constrained optimization problem and its continuous reformulation up to their type. For that, we focus on the nondegenerate M- and T-stationary points, respectively. Their so-called M- and T-indices, which uniquely determine the global and local structure of optimization problems under consideration in algebraic terms, are traced. As novelty, we suggest to regularize the continuous reformulation for this purpose. The main consequence of our analysis is that the number of saddle points of the regularized continuous reformulation grows exponentially as compared to that of the initial cardinality-constrained optimization problem. Additionally, we obtain the Morse theory for the regularized
continuous reformulation by using the corresponding results on mathematical programs with orthogonality type constraints.
}

\vspace{2ex}

\vspace{2ex}
{\bf MSC-classification: 90C26, 49M20}

\section{Introduction}

We consider the class of cardinality-constrained optimization problems:
\[
\mbox{CCOP}: \quad
\min_{x} \,\, f(x)\quad \mbox{s.\,t.} \quad h(x)=0, \quad g(x)\ge 0, \quad \left\|x\right\|_0\le s
\]
with the feasible set given by equality, inequality, and cardinality constraints, where the so-called zero ”norm” is counting non-zero entries of $x$:
\[
\left\|x\right\|_0 = \left|\left\{i \in \{1,\ldots, n\}\; \vert\; x_i\ne 0\right\}\right|.
\]
Here, we assume that the objective function $f \in C^2(\R^n,\R)$, as well as the equality and inequality constraints 
$h=\left(h_p, p \in P\right) \in C^2(\R^n,\R^{|P|})$,
$g=\left(g_q, q \in Q\right) \in C^2(\R^n,\R^{|Q|})$ are twice continuously differentiable, and $s \in \{0,1,\ldots, n-1\}$ is an integer. The bounded zero norm of the decision variable $x$ induces sparsity, which is motivated by various applications, such as compressed sensing, model selection, image processing etc., see e.\,g. \cite{donoho:2006}, \cite{tibshirani:1996}, and \cite{shechtman:2011}.

Let us start by recalling the continuous reformulation of CCOP from 
\cite{burdakov:2016}. There, auxiliary binary $y$-variables were used in order to rewrite the sparsity constraints. After relaxing the  binary constraints in a standard way, the authors arrive at
\begin{equation}
     \label{eq:relax}
     \begin{array}{rl}
    \displaystyle \min_{x,y} \,\, f(x) \quad \mbox{s.\,t.} \quad 
   &h(x)=0, \quad g(x)\ge 0, \\
   &\displaystyle \sum_{i=1}^{n} y_i \geq  n - s, \quad x_i y_i =0, \quad 0 \leq y_i \leq 1,  \quad 
         i=1, \ldots, n.
   \end{array}
\end{equation}
As pointed out in \cite{burdakov:2016}, $\bar x$ solves CCOP if and only if there exists a vector $\bar y$ such that $\left(\bar x, \bar y\right)$ solves (\ref{eq:relax}). 
The question on the global structure of the continuous reformulation (\ref{eq:relax}) and the initial CCOP arises quite naturally. For that, we focus on the stationary points -- minimizers, but also all kinds of saddle points -- which are topologically relevant in the sense of Morse theory, see e.g. \cite{Jongen:2000}. Morse theory typically provides deformation and cell-attachment results for the lower level sets. Deformation means that outside the set of stationary points the topology of lower level sets remains unchanged if the level varies. Cell-attachment describes the topological changes of the lower level sets if passing a nondegenerate stationary point.
As a consequence, it becomes possible to relate the numbers of minimizers and saddle points, thus, to adequately capture the global structure of the underlying optimization problem. While doing so, the notion of nedegeneracy becomes crucial. Nondegeneracy refers to some tailored versions of linear independence constraint qualification, strict complementarity and second-order regularity. It is justified by the fact that all stationary points are generically nondegenerate, i.e. this property holds on an open and dense subset of defining functions. Assuming nondegeneracy, stationary points can be classified according to their index which encodes the local structure of the optimization problem in algebraic terms. 

In \cite{laemmel:ccop}, M-stationary points were identified to adequately describe the global structure of CCOP. 
The attempt to establish Morse theory also for its continuous reformulation (\ref{eq:relax}) has been undertaken in our recent paper \cite{laemmel:mpoc}. There, the class of mathematical problems with orthogonality type constraints (MPOC) has been introduced for the latter purpose. By using the notion of T-stationarity, the deformation and cell-attachment results for MPOC were shown.
Although (\ref{eq:relax}) is a subclass of MPOC, and its T-stationary points naturally correspond to the M-stationary points of CCOP, the results from \cite{laemmel:mpoc} are not applicable. The reason is that all T-stationary points of (\ref{eq:relax}) turn out to be degenerate and, thus, the cell-attachment cannot be performed. This is even the case if we start from a CCOP with all its M-stationary points being nondegenerate.
In order to nevertheless examine the global structure of (\ref{eq:relax}), we suggest to study its regularized version. Our main idea is not only to linearly perturb the objective function in (\ref{eq:relax}) with respect to $y$-variables, but also to additionally relax the upper bounds on them:
\[
     \begin{array}{rl}
    \displaystyle \mathcal{R}(c,\varepsilon): \quad \min_{x,y} \,\, f(x) +c^Ty\quad \mbox{s.\,t.} 
   &h(x)=0, \quad g(x)\ge 0, \\
   &\displaystyle \sum_{i=1}^{n} y_i \geq  n - s, \quad x_i y_i =0, \quad 0 \leq y_i \leq 1+\varepsilon,  \quad 
         i=1, \ldots, n,
   \end{array}
\]
where $c \in \R^n$ and $\varepsilon>0$.
By doing so, the proposed regularization $\mathcal{R}$ remains an MPOC, whose T-stationary points now become generically nondegenerate. This follows from the derived relation between M-stationary points of CCOP and T-stationary points of $\mathcal{R}$. We also successfully trace their corresponding M- and T-indices by imposing -- except of nondegeneracy -- an additional condition on $\mathcal{R}$. The latter requires that in presence of biactive orthogonality type constraints, the multipliers of a T-stationary point corresponding to zero $x$-variables and nonzero $y$-variables do not vanish. Fortunately, this additional condition can be shown to hold generically, thus, it is far from being restrictive.  
As a consequence of our analysis, the Morse theory for the regularized continuous reformulation $\mathcal{R}$ can be deduced from the corresponding results on MPOC from \cite{laemmel:mpoc}. More importantly, we estimate the
number of its saddle points. Namely, each saddle point of CCOP generates exponentially many saddle points of $\mathcal{R}$, all of them having the same index. 
 
The article is organized as follows. In Section 2 we discuss some preliminary notions for CCOP. 
Section 3 is devoted to the analysis of the regularized continuous reformulation $\mathcal{R}$.
In Section 4 we compare our findings with those for the continuous reformulation (\ref{eq:relax}) known from the literature.

Our notation is standard. The cardinality of a finite set $A$ is denoted by $|A|$. 
The $n$-dimensional
Euclidean space is denoted by $\mathbb{R}^n$ with the coordinate vectors $e_i,i= 1,\ldots, n$. 
Given a twice continuously differentiable function $f:\mathbb{R}^n\rightarrow \mathbb{R}$, $\nabla f$ denotes its gradient, and $D^2f$ stands for its Hessian. 

\section{Cardinality-constrained optimization problems}
\label{sec:ccop}

We start by recalling details on the class of CCOP from the literature. The following notation for a CCOP feasible point $\bar x$ will be used:
\[
Q_0(\bar x)=\left\{q \in Q \, \left\vert\, g_q(\bar x)=0\right.\right\}
\]
denotes the index set of active inequality constraints, and the index set of its vanishing components is set to be 
\[
I_0(\bar x)=\left\{i \in \{1,\ldots, n\}\, \left\vert\, \bar x_i= 0\right.\right\}.
\]

Let us mention the CCOP-tailored linear independence constraint qualification. It is known to hold generically on the whole CCOP feasible set, see \cite{laemmel:ccop}.

\begin{definition}[CC-LICQ, see \cite{cervinka:2016}]
\label{def:CC-LICQ}
We say that a CCOP feasible point $\bar x$ of CCOP satisfies the cardinality-constrained linear independence constraint qualification (CC-LICQ)  if the following gradients are linearly independent:
\[
\nabla h_p(\bar x), p \in P, \quad
\nabla g_q(\bar x), q \in Q_0(\bar x), \quad
e_i, i \in I_0(\bar x).
\]
\end{definition}

The topologically relevant concept of M-stationarity for CCOP is stated as follows.

\begin{definition}[M-stationarity, see \cite{burdakov:2016}]
\label{def:M-stat}
A CCOP feasible point $\bar x$ is called M-stationary if there exist multipliers
  $$\bar \lambda_p, p \in P, \bar \mu_q, q \in Q_0(\bar x), \bar \gamma_i, i \in I_0(\bar x),$$
such that the following conditions hold:
\begin{equation}
\label{eq:mstat-1}
\nabla f(\bar x) = \sum\limits_{p \in P}\bar \lambda_p \nabla h_p(\bar x)+
    \sum\limits_{q \in Q_0(\bar x)}\bar \mu_q \nabla g_q(\bar x)
    +\sum\limits_{i\in I_0(\bar x)} \bar \gamma_i e_i,
\end{equation}
\begin{equation}
\label{eq:mstat-2}
\bar \mu_q \geq 0 \mbox{ for all }q \in Q_0(\bar x).
\end{equation}
\end{definition}
Since the multipliers are unique under CC-LICQ, it is convenient to define the Lagrange function:
    \[
     L(x)=f\left(x\right) - \sum\limits_{p \in P}\bar \lambda_p h_p\left(x\right)-
    \sum\limits_{q \in Q_0(\bar x)}\bar \mu_q g_q\left(x\right)
    -\sum\limits_{i \in I_0(\bar x)} \bar \gamma_i x_i.
    \]
We also use the corresponding tangent space:
\[
    \mathcal{T}_{\bar x}=\left\{
\xi \in \R^n\,\left\vert\, Dh_p(\bar x) \xi=0, p \in P, 
Dg_q(\bar x)\xi=0,q \in Q_0(\bar x),
\xi_i=0, i \in I_0(\bar x)
\right.\right\}.
\]
Let us proceed with the definition of nondegeneracy for M-stationary points as introduced in \cite{laemmel:ccop}. It is justified there by showing that all M-stationary points of CCOP are generically nondegenerate.
    
\begin{definition}[Nondegenerate M-stationarity, see \cite{laemmel:ccop}]
An M\Hyphdash stationary point $\bar x$ of CCOP is called nondegenerate if
\begin{itemize}
    \item[] NDM1: CC-LICQ holds at $\bar x$,
    \item[] NDM2: $\bar \mu_q>0$ for all $q\in Q_0(\bar x)$,
    \item[] NDM3: if $\left\|\bar x\right\|_0<s$ then $\bar \gamma_i\ne 0$ for all $i\in I_0(\bar x)$,
    \item[] NDM4: the matrix $D^2 L(\bar x)\restriction_{\mathcal{T}_{\bar x}}$ is nonsingular.  
    \end{itemize}
\end{definition}
The nondegeneracy conditions NDM1-NDM4 are tailored for the CCOP class. Note that NDM2 corresponds to the strict complementarity and NDM4 to the second-order regularity as they are typically defined in the context of nonlinear programming. NDM1 substitutes the usual linear independence constraint qualification. NDM3 is new and says that, unless the sparsity constraint is active, the corresponding Lagrange multipliers should not vanish.
With a nondegenerate M-stationary point $\bar x$ an M-index can be associated. M-index captures the structure of CCOP locally around $\bar x$ and defines the type of an M-stationary point, see \cite{laemmel:ccop} for details. In particular, nondegenerate minimizers of CCOP are characterized by the vanishing M-index. If M-index does not vanish, we get all kinds of saddle points.  

\begin{definition}[M-Index, see \cite{laemmel:ccop}]
\label{def:mindex}
Let $\bar x$ be a nondegenerate M\Hyphdash stationary point of CCOP.
The  number of negative eigenvalues of the matrix $D^2 L(\bar x)\restriction_{\mathcal{T}_{\bar x}}$ is called its quadratic index ($QI$). 
The number $s-\left\|\bar x\right\|_0$ is called the sparsity index ($SI$) of $\bar x$. We define the M-index ($MI$) as the sum of both, i. e. $MI=SI+BI$.
\end{definition}

Let us mention the role M-stationary points and their M-indices play for capturing the global structure of CCOP. In \cite{laemmel:ccop}, deformation and cell-attachment in the sense of Morse theory are proved for a generic CCOP. Deformation says that lower level sets are homeomorphic if passing a level which does not correspond to any M-stationary point. Cell-attachment algebraically describes topological differences between lower level sets if a level corresponding to a nondegenerate M-stationary is crossed. Exponentially many cells of the same dimension need to be attached to a lower level set in order to obtain another lower level set up to homotopy-equivalence. The dimension of those cells to be attached coincides with the M-index of a nondegenerate M-stationary point. A global interpretation of deformation and cell-attachment can be given in form of a mounting pass result, see again \cite{laemmel:ccop} for details.

\section{Regularized continuous reformulation}

Let us associate with CCOP the regularized continuous reformulation
, cf. (\ref{eq:relax}):
\[
     \begin{array}{rl}
    \displaystyle \mathcal{R}(c,\varepsilon): \quad \min_{x,y} \,\, f(x) +c^Ty\quad \mbox{s.\,t.} 
   &h(x)=0, \quad g(x)\ge 0, \\
   &\displaystyle \sum_{i=1}^{n} y_i \geq  n - s, \quad x_i y_i =0, \quad 0 \leq y_i \leq 1+\varepsilon,  \quad 
         i=1, \ldots, n,
   \end{array}
\]
where $c \in \R^n$ and $\varepsilon>0$.
Given a feasible point  $(\bar x,\bar y)$ of $\mathcal{R}$, we define the index sets which correspond to the orthogonality type constraints $x_i y_i=0$, $y_i \geq 0$, $i=1, \ldots, n$:
    \[
    a_{00}\left(\bar x,\bar y\right)=\left\{i\in\left\{1,\ldots,n\right\}\,\left\vert\, \bar x_i=0,\bar y_i=0\right.\right\},
    \]    
         \[
     a_{01}\left(\bar x,\bar y\right)=\left\{i\in\left\{1,\ldots,n\right\}\,\left\vert\, \bar x_i=0,\bar y_i>0\right.\right\},
     \]   
    \[
    a_{10}\left(\bar x,\bar y\right)=\left\{i\in\left\{1,\ldots,n\right\}\,\left\vert\, \bar x_i\ne 0,\bar y_i=0\right.\right\}.
    \]
The $y$-components which attain the upper bound are stored in  
     \[\mathcal{E}(\bar y)=\left\{i \in\{1,\ldots,n\}\left\vert\,\bar y_i=1+\varepsilon\right.\right\}.\]
The index set of the active inequality constraints remains to be denoted by
\[
Q_0(\bar x)=\left\{q \in Q \, \left\vert\, g_q(\bar x)=0\right.\right\}.
\]

The following assumption on the regularization parameters $c$ and $\varepsilon$ will be helpful in what follows. Note that it is not restrictive since a randomly generated vector $c$ with positive components fulfills Assumption \ref{ass:c} with probability one. The upper bound on $\varepsilon$ depends just on $n$ and $s$, but not on the CCOP defining functions. If not stated otherwise, Assumption \ref{ass:c} holds throughout the whole paper.

\begin{assumption}
\label{ass:c}
Let the components of $c$ be positive and pairwise different, and $\varepsilon \leq \frac{1}{n-s}$.
\end{assumption}

The regularized continuous reformulation $\mathcal{R}$ is a special case of mathematical programs with orthogonality type constraints (MPOC). The latter class was examined in \cite{laemmel:mpoc}, where the MPOC-tailored linear independence constraint qualification, the topologically relevant notion of (nondegenerate) T-stationary points with the corresponding T-index were introduced. 
We apply these concepts to the regularization $\mathcal{R}$. 

\begin{definition}[MPOC-LICQ]
We say that a feasible point $(\bar x,\bar y)$ of $\mathcal{R}$ satisfies the MPOC-tailored linear independence constraint qualification (MPOC-LICQ) if the following vectors are linearly independent:
\[\begin{pmatrix}
\nabla h_p(\bar x)\\0
\end{pmatrix}, p \in P, \quad 
\begin{pmatrix}
\nabla g_q(\bar x)\\0
\end{pmatrix}, q \in Q_0(\bar x), \quad  
\begin{pmatrix}
0\\
e_i
\end{pmatrix},i\in \mathcal{E}(\bar y), \quad
\begin{pmatrix}
0\\
e
\end{pmatrix}
\mbox{ if } \sum\limits_{i=1}^{n} \bar y_i = n - s,
\]
\[
\begin{pmatrix}
e_i\\
0
\end{pmatrix}, i \in a_{01}\left(\bar x,\bar y\right)\cup a_{00}\left(\bar x,\bar y\right), \quad 
\begin{pmatrix}
0\\
e_i
\end{pmatrix}, i \in a_{10}\left(\bar x, \bar y\right)\cup a_{00}\left(\bar x, \bar y\right).
\]
\end{definition}

Let us relate CCOP- and MPOC-tailored linear independence constraint qualifications. 

\begin{theorem}[CC-LICQ vs. MPOC-LICQ]
\label{thm:cc-mpoc-licq}
A feasible point $\bar x$ of CCOP fulfills
CC-LICQ if and only if MPOC-LICQ holds at any feasible point $\left(\bar x, \bar y\right)$ of
$\mathcal{R}$.
\end{theorem}
\proof
CC-LICQ holds at $\bar x$ if the following vectors are linearly independent:
\[
\nabla h_p(\bar x), p \in P, \quad
\nabla g_q(\bar x), q \in Q_0(\bar x), \quad
e_i, i \in I_0(\bar x).
\]
Since $I_0(\bar x)=a_{01}\left(\bar x,\bar y\right)\cup a_{00}\left(\bar x,\bar y\right)$, the vectors
\[\begin{pmatrix}
\nabla h_p(\bar x)\\0
\end{pmatrix}, p \in P,\quad
\begin{pmatrix}
\nabla g_q(\bar x)\\0
\end{pmatrix}, q \in Q_0(\bar x),\quad
\begin{pmatrix}
e_i\\
0
\end{pmatrix}, i \in a_{01}\left(\bar x,\bar y\right)\cup a_{00}\left(\bar x,\bar y\right),
\]
are linearly independent.
In order to prove that MPOC-LICQ is indeed fulfilled at $\left(\bar x, \bar y\right)$, it remains therefore to show that the following vectors are also linearly independent:
\begin{equation}
\label{eq:vectors}
\begin{pmatrix}
0\\
e_i
\end{pmatrix},i\in \mathcal{E}(\bar y),\quad 
\begin{pmatrix}
0\\
e
\end{pmatrix}\mbox{ if } \sum\limits_{i=1}^{n} \bar y_i =n - s, \quad
\begin{pmatrix}
0\\
e_i
\end{pmatrix}, i \in a_{10}\left(\bar x, \bar y\right)\cup a_{00}\left(\bar x, \bar y\right).
\end{equation}
Since $\mathcal{E}\left(\bar y\right)\subset a_{01}\left(\bar x,\bar y\right)$, we have that
$\begin{pmatrix}
0\\
e_{i_1}
\end{pmatrix}$ and 
$\begin{pmatrix}
0\\
e_{i_2}
\end{pmatrix}$
are linearly independent for all $i_1 \in \mathcal{E}\left(\bar y\right)$ and all
$i_2 \in a_{10}\left(\bar x, \bar y\right)\cup a_{00}\left(\bar x, \bar y\right)$.
It is enough to ensure that we have at most $n$ vectors in (\ref{eq:vectors}). We distinguish the following cases:

a) $\sum\limits_{i=1}^{n} \bar y_i >n - s$. Due to the feasibility of $(\bar x, \bar y)$, we have
\[
\left|a_{01}\left(\bar x, \bar y\right)\right|+
\left|a_{10}\left(\bar x, \bar y\right)\right|+
\left|a_{00}\left(\bar x, \bar y\right)\right|=n.
\]
Hence, MPOC-LICQ holds if and only if
$\left|\mathcal{E}\left(\bar y\right)\right|\le\left|a_{01}\left(\bar x, \bar y\right)\right|$,
which is obviously true.

b) $\sum\limits_{i=1}^{n} \bar y_i =n - s$.
Assume that MPOC-LICQ does not hold at $(\bar x, \bar y)$. Consequently, we must have $\mathcal{E}\left(\bar y\right)=a_{01}\left(\bar x, \bar y\right)$.
However, we obtain then:
\[
n-s=\sum\limits_{i=1}^{n} \bar y_i=\left|a_{01}\left(\bar x, \bar y\right)\right|\left(1+\varepsilon\right)=\left|a_{01}\left(\bar x, \bar y\right)\right|+
\left|a_{01}\left(\bar x, \bar y\right)\right|\varepsilon.
\]
In particular, $\left|a_{01}\left(\bar x, \bar y\right)\right|\varepsilon$ is an integer.
Due to $\varepsilon \le \frac{1}{n-s}$ it must, thus, hold either
$\left|a_{01}\left(\bar x, \bar y\right)\right|=0$ or $\left|a_{01}\left(\bar x, \bar y\right)\right| \ge n-s$.
The former implies 
\[
\sum\limits_{i=1}^{n} \bar y_i=0< n-s,
\]
whereas the latter implies
\[
\sum\limits_{i=1}^{n} \bar y_i\ge (n-s)(1+\varepsilon)>n-s.
\]
Both cases yield a contradiction. The reverse implication of the assertion is straightforward.
\qed

Now, we turn our attention to the topologically relevant concept of T-stationarity  for $\mathcal{R}$.

\begin{definition}[T-stationary point]
\label{def:t-stat}
A feasible point $(\bar x,\bar y)$ of $\mathcal{R}$ is called T-stationary if there exists multipliers 
\[
\begin{array}{l}
\bar \lambda_p, p \in P,
\bar \mu_{1,q}, q\in Q_0(\bar x),\bar \mu_{2,i},i \in \mathcal{E}(\bar y),\bar \mu_3,\\ \bar \sigma_{1,i_{01}}, i_{01} \in a_{01}\left(\bar x,\bar y\right),\bar \sigma_{2,i_{10}}, i_{10} \in a_{10}\left(\bar x,\bar y\right),\bar \varrho_{1,i_{00}},\bar \varrho_{2,i_{00}}, i_{00} \in a_{00}\left(\bar x, \bar y\right),
\end{array}
\]
such that the following conditions hold:
\begin{equation}
   \label{eq:tstat-1} 
   \begin{array}{rcl}
  \begin{pmatrix}
   \nabla f(\bar x)\\
   c
  \end{pmatrix}&=& \displaystyle\sum\limits_{p \in P}\bar \lambda_p 
  \begin{pmatrix}
  \nabla h_p(\bar x)\\
  0
  \end{pmatrix}+
    \sum\limits_{q \in Q_0(\bar x)}\bar \mu_{1,q} 
    \begin{pmatrix}
    \nabla g_q(\bar x)\\
    0
    \end{pmatrix}-
    \sum\limits_{i\in\mathcal{E}(\bar y)} \bar \mu_{2,i} \begin{pmatrix}
    0\\
    e_i
\end{pmatrix}+
\bar \mu_3 \begin{pmatrix}
0\\
e
\end{pmatrix}\\ \\
    && \displaystyle +\sum\limits_{i_{01} \in a_{01}\left(\bar x,\bar y\right)} \bar \sigma_{1,i_{01}}
    \begin{pmatrix}
    e_{ i_{01}}\\
    0
    \end{pmatrix}
    +\sum\limits_{i_{10} \in a_{10}\left(\bar x,\bar y \right)} \bar \sigma_{2,i_{10}}    \begin{pmatrix}
    0\\
    e_{ i_{10}}
    \end{pmatrix}\\ \\
    &&\displaystyle+\sum\limits_{i_{00} \in a_{00}\left(\bar x, \bar y\right)} \left(\bar \varrho_{1,i_{00}}
    \begin{pmatrix}
   e_{ i_{00}}\\
    0
    \end{pmatrix}
   +\bar \varrho_{2,i_{00
    }}
    \begin{pmatrix}
    0\\
    e_{ i_{00}}
    \end{pmatrix}\right), \end{array}
\end{equation}
\begin{equation}
   \label{eq:tstat-2} \bar \mu_{1,q} \ge 0 \mbox{ for all } q\in Q_0\left(\bar x\right),
   \bar \mu_{2,i} \ge 0 \mbox{ for all } i\in \mathcal{E}(\bar y),\bar \mu_3\ge 0, \bar\mu_3\cdot\left(
   \sum\limits_{i=1}^{n} \bar y_i -(n-s)\right)=0,
\end{equation}
\begin{equation}
   \label{eq:tstat-3} \bar \varrho_{1,i_{00}}=0 \mbox{ or }\bar \varrho_{2,i_{00}}\le 0 \mbox{ for all } i_{00} \in a_{00}\left(\bar x,\bar y\right).
\end{equation}
\end{definition}

We again define the appropriate Lagrange function:
\[
        \begin{array}{rcl}
        L^\mathcal{R}(x,y)&=& \displaystyle
   f(x)+
   c^T y -
        \sum\limits_{p \in P}\bar \lambda_p
  h_p( x)-
    \sum\limits_{q \in Q_0(\bar x)}\bar \mu_{1,q} 
    g_q(x) +\sum\limits_{i\in\mathcal{E}(\bar y)} \bar \mu_{2,i} 
   \left( y_i- (1+\varepsilon)\right)\\ \\ & &\displaystyle
-
\bar \mu_3 
 \left(\sum\limits_{i=1}^{n} y_i - (n-s)\right)
      -\sum\limits_{i_{01} \in a_{01}\left(\bar x, \bar y\right)} \bar \sigma_{1,i_{01}}
     x_{ i_{01}}
    -\sum\limits_{i_{10} \in a_{10}\left(\bar x, \bar y\right)} \bar \sigma_{2,i_{10}}   
     y_{ i_{10}}
   \\ \\
    &&\displaystyle -\sum\limits_{i_{00} \in a_{00}\left(\bar x, \bar y \right)} \left(\bar \varrho_{1,i_{00}}
    x_{ i_{00}}
   +\bar \varrho_{2,i_{00
    }}
     y_{ i_{00}}
    \right).
    \end{array}
\]
Moreover, we set for the corresponding tangential space:
\[
    \mathcal{T}^{\mathcal{R}}_{(\bar x,\bar y)} =\left\{
\xi \in \R^{2n}\,\left\vert\, \begin{array}{l} \begin{pmatrix}

Dh_p(\bar x), 0\end{pmatrix} \xi=0, p \in P, 
\begin{pmatrix}
Dg_q(\bar x), 0\end{pmatrix}\xi=0,q \in Q_0(\bar x),\\
\begin{pmatrix}
0,e_i
\end{pmatrix}\xi=0, i\in \mathcal{E}(\bar y), 
\begin{pmatrix}
0,e
\end{pmatrix}\xi=0 \mbox{ if } \displaystyle \sum_{i=1}^{n} \bar y_i = n - s,
\\
\begin{pmatrix}
e_i,0
\end{pmatrix}\xi=0, i \in a_{00}(\bar x,\bar y) \cup a_{01}(\bar x,\bar y),\\
\begin{pmatrix}
0,e_i
\end{pmatrix}\xi=0, i \in a_{00}(\bar x,\bar y) \cup a_{10}(\bar x,\bar y)
\end{array}
\right.\right\}.
\]

\begin{definition}[Nondegenerate T-stationary point]
A T-stationary point $(\bar x,\bar y)$ of $\mathcal{R}$ with multipliers $(\bar \lambda, \bar \mu, \bar \sigma, \bar \varrho)$ is called nondegenerate if
\begin{itemize}
    \item []NDT1: MPOC-LICQ holds at $(\bar x,\bar y)$,
    \item []NDT2: the strict complementarity (SC) holds for active inequality constraints, i.\,e. $\bar \mu_{1,q}>0$ for all $q \in Q_0\left(\bar x\right)$, $\bar \mu_{2,i}>0$ for all $i \in \mathcal{E}\left(\bar y\right)$, and if  $\sum\limits_{i=1}^{n} \bar  y_i =n-s$ then also $\bar \mu_3>0$,
    \item []NDT3: the multipliers corresponding to biactive orthogonality type constraints do not vanish, i.\,e. $\bar\varrho_{1,i_{00}}\ne 0$ and $\bar\varrho_{2,i_{00}}< 0$ for all $i_{00}\in a_{00}\left(\bar x,\bar y\right)$,
    \item []NDT4: the matrix $D^2 L^{\mathcal{R}}(\bar x,\bar y)\restriction_{\mathcal{T}^\mathcal{R}_{(\bar x,\bar y)}}$ is nonsingular.
\end{itemize}  
For a nondegenerate T-stationary point we eventually use an additional condition:
\begin{itemize}
\item [] NDT5: if $a_{00}\left(\bar x,\bar y\right)\not = \emptyset$, then $\bar \sigma_{1,i_{01}}\ne 0$ for all $i_{01}\in a_{01}(\bar x, \bar y)$.
\end{itemize}

\end{definition}

\begin{definition}[T-index]
Let $(\bar x,\bar y)$ be a nondegenerate T-stationary point of $\mathcal{R}$ with unique multipliers $\left(\bar \lambda,\bar \mu, \bar \sigma,\bar \varrho\right)$. The number of negative eigenvalues of the matrix $D^2 L^{\mathcal{R}}(\bar x,\bar y)\restriction_{\mathcal{T}^\mathcal{R}_{(\bar x,\bar y)}}$ is called its quadratic index ($QI$). The cardinality of $a_{00}\left(\bar x,\bar y\right)$ is called the biactive index ($BI$) of $(\bar x,\bar y)$. We define the T-index ($TI$) as the sum of both, i.\,e. $TI=QI+BI$.
\end{definition}

The following Lemma \ref{lem:a01} provides insights into the structure of auxiliary $y$-variables corresponding to a T-stationary point of $\mathcal{R}$.

\begin{lemma}[Auxiliary $y$-variables in $\mathcal{R}$]
\label{lem:a01} 
Let $(\bar x,\bar y)$ be a T-stationary point of $\mathcal{R}$, then it holds:
\begin{itemize}
    \item [a)] the summation inequality constraint is active, i.\,e. $ \sum\limits_{i=1}^{n} \bar y_i =n - s$,
    \item [b)] the index set $a_{01}(\bar x,\bar y)$ consists of exactly $n-s$ elements,
    \item [c)] $n-s-1$ components of $\bar y$ are equal to $1+\varepsilon$, one component is equal to
    $1-(n-s-1)\varepsilon$, and $s$ remaining components vanish.
\end{itemize}
\end{lemma}
\proof
 a)
Let $(\bar x,\bar y)$ be a T-stationary point of $\mathcal{R}$
and $ \sum_{i=1}^{n} \bar y_i >n - s$. Then, there exist multipliers
$(\bar \lambda, \bar \mu, \bar \sigma, \bar \varrho)$
such that (\ref{eq:tstat-1})--(\ref{eq:tstat-3}) are fulfilled.
Since $\bar \mu_3=0$, we have that the $(n+i)$-th row of (\ref{eq:tstat-1}) reads as
\[
c_i=\left\{
\begin{array}{ll}
     -\bar \mu_{2,i},&\mbox{for }i \in \mathcal{E}(\bar y),  \\
     \bar \sigma_{2,i}, &\mbox{for }i \in a_{10}\left(\bar x,\bar y\right),\\
     \bar \varrho_{2,i}, &\mbox{for }i\in a_{00}\left(\bar x,\bar y\right),\\
     0,&\mbox{else}.
\end{array}\right.
\]
Due to $c>0$, the sets $a_{01}\left(\bar x,\bar y\right)$ and $\mathcal{E}(\bar y)$ have to be equal and, moreover,
due to $\bar \mu_{2,i}\ge 0$, they have to be empty. But then, clearly, $ \sum_{i=1}^{n} \bar  y_i =0<n - s$, a contradiction.

b)
Since $(\bar x,\bar y)$ is a T-stationary point, there exist
$(\bar \lambda, \bar \mu, \bar \sigma, \bar \varrho)$
such that (\ref{eq:tstat-1})--(\ref{eq:tstat-3}). By the proof of statement a) we can conclude that
$\bar \mu_3 >0$. Hence,
the $(n+i)$-th row reads as
\begin{equation}
    \label{eq:k-throw}
c_i=\left\{
\begin{array}{ll}
     -\bar \mu_{2,i}+\bar \mu_3,&\mbox{for }i \in \mathcal{E}(\bar y),  \\
     \bar \sigma_{2,i}+\bar \mu_3, &\mbox{for }i \in a_{10}\left(\bar x,\bar y\right),\\
     \bar \varrho_{2,i}+\bar \mu_3, &\mbox{for }i\in a_{00}\left(\bar x,\bar y\right),\\
     \bar \mu_3,&\mbox{else}.
\end{array}\right.
\end{equation}
Let us assume that the index set $a_{01}(\bar x,\bar y)$ consists of fewer than $n-s$ elements. Then, we have by using Assumption \ref{ass:c}:
\[ \sum_{i=1}^{n} \bar y_i \le (n-s-1)\cdot (1+\varepsilon)\le n-s-1+\frac{n-s-1}{n-s} <n - s,
\]
a contradiction to feasibility. 
Let us assume that the index set $a_{01}(\bar x,\bar y)$ consists of more than $n-s$ elements instead.
Since the components of $c$ are assumed to be pairwise different, we see from (\ref{eq:k-throw}) that there exists at most one element in $a_{01}\left(\bar x,\bar y\right)\backslash \mathcal{E}(\bar y)$ and, consequently, there are at least $n-s$ elements in $\mathcal{E}(\bar y)$. Therefore, we have:
\[
\sum_{i=1}^{n} \bar y_i \geq (n-s)\cdot (1+\varepsilon) > n-s,
\]
which contradicts a). 

c)
Due to b), $a_{01}(\bar x,\bar y)$ consists of exactly $n-s$ elements. We conclude as in b) that there is at most one element in $a_{01}\left(\bar x,\bar y\right)\backslash \mathcal{E}(\bar y)$. In view of statement a), $\mathcal{E}(\bar y)$ cannot consist of $n-s$ elements and, thus, must consist of $n-s-1$ elements. Hence, the statement follows immediately.
\qed

We are ready to identify how many T-stationary points of $\mathcal{R}$ are generated by an M-stationary point of CCOP, and of what type they are.


\begin{theorem}[Stationarity from CCOP to $\mathcal{R}$]
\label{thm:mtot}
 If $\bar x$ is an M-stationary point of CCOP, then there exist at least $\binom{n-\left\|\bar x\right\|_0-1}{n-s-1}$ choices of $\bar y$ such that $(\bar x,\bar y)$ is a T-stationary point of $\mathcal{R}$. 
If $\bar x$ is additionally nondegenerate with M-index $m$, then all corresponding T-stationary points $(\bar x, \bar y)$ are also nondegenerate with T-index $m$. Moreover, their number is exactly $\binom{n-\left\|\bar x\right\|_0-1}{n-s-1}$, and NDT5 holds at any of them.
\end{theorem}
\proof
Since $\bar x$ is an M-stationary point, there exist multipliers
$(\bar \lambda, \bar \mu,\bar \gamma)$ with
(\ref{eq:mstat-1}), (\ref{eq:mstat-2}). We set $$\bar i=\mbox{argmax}\left\{c_i\left\vert i\in I_0\left(\bar x\right)\right.\right\}$$ to be the index of $I_0\left(\bar x\right)$ for which
$c_i$ is maximal. Furthermore, let $ \mathcal{\bar E}\subset{I_0}\left(\bar x\right)\backslash\left\{\bar i\right\}$ be an index subset with $n-s-1$ elements, i.\,e. $\left\vert\mathcal{\bar E}\right\vert = n-s-1$. Note that this is always possible since $I_0\left(\bar x\right)$ consists of at least $n-s$ elements.  Next, we set
\[
\bar y_i=\left\{
\begin{array}{ll}
1+\varepsilon,&\mbox{for }i \in \mathcal{\bar E},\\
1-(n-s-1)\varepsilon,&\mbox{for }i=\bar i,\\
0,&\mbox{else}.
\end{array}\right.
\] 
Consequently, the point $(\bar x,\bar y)$ is feasible for $\mathcal{R}$.
Let the multipliers $\bar \lambda_p$, $p \in P$, corresponding to the equality constraints remain unchanged. We set the other multipliers as follows: 
\[
\bar \mu_{1,q}=\bar \mu_q, q \in Q_0(\bar x),\quad \bar \mu_{2,i}=c_{\bar i}-c_i,i \in \mathcal{\bar E},\quad
\bar \mu_3=c_{\bar i},
\]
\[
\bar \sigma_{1,i_{01}}=\bar\gamma_{i_{01}},i_{01}\in  \mathcal{\bar E} \cup \left\{\bar i\right\},\quad
\bar \sigma_{2,i_{10}}=c_{i_{10}}-c_{\bar i},i_{10}\in \{1,\ldots,n\}\backslash I_0\left(\bar x\right),
\] 
\[
\bar \varrho_{1,i_{00}}=\bar\gamma_{i_{00}}, \quad \bar \varrho_{2,i_{00}}=c_{i_{00}}-c_{\bar i}, i_{00}\in I_0\left(\bar x\right)\backslash\left(\mathcal{\bar E}\cup \left\{\bar i\right\}\right).
\]
We note that 
\[
\mathcal{E}(\bar y)= \mathcal{\bar E}, \quad a_{01}(\bar x,\bar y)=\mathcal{\bar E}\cup \left\{\bar i\right\}, \quad
a_{10}(\bar x,\bar y)=\{1,\ldots,n\}\backslash I_0\left(\bar x\right), \quad
a_{00}(\bar x,\bar y)=I_0\left(\bar x\right)\backslash\left(\mathcal{\bar E}\cup \left\{\bar i\right\}\right).
\]
Since $\bar x$ is M-stationary and  $\sum\limits_{i=1}^{n} \bar y_i =n - s$, we obtain the T-stationarity condition (\ref{eq:tstat-1}) for $(\bar x, \bar y)$.
Let us check the signs of the multipliers. Due to the M-stationarity of $\bar x$ or by construction, we have:
\[  
\bar \mu_{1,q} \geq 0, q \in Q_0(\bar x), \quad \bar \mu_{2,i} \geq 0,i\in \mathcal{\bar E}, \quad \bar \mu_3 \geq 0, \quad \bar \varrho_{2,i_{00}}\leq 0, i_{00} \in I_0\left(\bar x\right)\backslash\left(\mathcal{\bar E}\cup \left\{\bar i\right\}\right).
\]
Thus, T-stationarity conditions (\ref{eq:tstat-2}), (\ref{eq:tstat-3}) are also fulfilled. 

It remains to show that there are at least $\binom{n-\left\|x\right\|_0-1}{n-s-1}$ possibilities to choose $\bar y$. This is exactly the number of possible choices for $\mathcal{\bar E}$. Hence, the number of possible choices for $\bar y$ cannot be less.
We assume that for $\bar x$ nondegenerate there is another $\tilde y$ which cannot be constructed as above, but so that $(\bar x,\tilde y)$ is a T-stationary point of $\mathcal{R}$. Obviously, we have $\tilde y_i=0$ for $i \notin I_0(\bar x)$. Additionally, we know due to Lemma \ref{lem:a01}a) that 
$\sum\limits_{i=1}^{n} \tilde y_i =n - s$ as well as from Lemma \ref{lem:a01}c) that there exists exactly one $\tilde i$ such that $0<\tilde y_{\tilde i}<1+\varepsilon$. Since $(\bar x,\tilde y)$ is T-stationary, conditions (\ref{eq:tstat-1})--(\ref{eq:tstat-3}) hold with multipliers $(\tilde \lambda, \tilde \mu, \tilde \sigma, \tilde \varrho)$. According to the aforementioned, the $(n+i)$-th row of (\ref{eq:tstat-1}) reads as
\[
c_i=\left\{
\begin{array}{ll}
     -\tilde \mu_{2,i}+\tilde \mu_3,&\mbox{for }i \in \mathcal{E}(\tilde y),  \\
       \tilde \mu_3,&\mbox{for }i=\tilde i,\\
     \tilde \sigma_{2,i}+\tilde \mu_3, &\mbox{for }i \in a_{10}\left(\bar x,\tilde y\right),\\
     \tilde \varrho_{2,i}+\tilde \mu_3, &\mbox{for }i\in a_{00}\left(\bar x,\tilde y\right).\\
\end{array}\right.
\]
Due to (\ref{eq:tstat-2}) we have that $c_i\le c_{\tilde i}$ for all $i \in \mathcal{E}(\tilde y)$. 
We show that it also holds $c_i\le c_{\tilde i}$ for all $i \in a_{00}\left(\bar x,\tilde y\right)$. 
In case of $\left\|\bar x\right\|_0=s$ we have $\left|I_0(\bar x)\right|=n-s=|a_{01}\left(\bar x,\tilde y\right)|$, where the second equality follows from Lemma \ref{lem:a01}b). Consequently, it follows from 
$I_0(\bar x)=a_{01}\left(\bar x,\tilde y\right)\cup a_{00}\left(\bar x,\tilde y\right)$
that $a_{00}(\bar x, \tilde y)=\emptyset$. 
Instead we suppose $\left\|\bar x\right\|_0<s$. Since $\bar x$ is a nondegenerate M-stationary point we have
\[
\nabla f(\bar x) = \sum\limits_{p \in P}\bar \lambda_p \nabla h_p(\bar x)+
    \sum\limits_{q \in Q_0(\bar x)}\bar \mu_q \nabla g_q(\bar x)
    +\sum\limits_{i_{01}\in a_{01}\left(\bar x,\tilde y\right)} \bar \gamma_{i_{01}} e_{i_{01}}
        +\sum\limits_{i_{00}\in a_{00}\left(\bar x,\tilde y\right)} \bar \gamma_{i_{00}} e_{i_{00}},
\]
where $\bar \gamma_i\ne 0$ for all $i\in I_0(\bar x)=a_{01}\left(\bar x, \tilde y\right)\cup a_{00}\left(\bar x,\tilde y\right).$
Additionally, we have due to $\left(\bar x,\tilde y\right)$ being T-stationary:
\[
\nabla f(\bar x) = \sum\limits_{p \in P}\tilde \lambda_p \nabla h_p(\bar x)+
    \sum\limits_{q \in Q_0(\bar x)}\tilde \mu_{1,q} \nabla g_q(\bar x)
    +\sum\limits_{i_{01}\in a_{01}\left(\bar x,\tilde y\right)} \tilde \sigma_{1,i_{01}} e_{i_{01}}+     \sum\limits_{i_{00}\in a_{00}\left(\bar x,\tilde y\right)} \tilde \varrho_{1,i_{00}} e_{i_{00}}.
\]
Then, CC-LICQ implies $\tilde \varrho_{1,i_{00}}=\bar \gamma_{i_{00}}\ne 0$ for all
$i\in a_{00}\left(\bar x,\tilde y\right)$. But then (\ref{eq:tstat-3}) implies
$\tilde \varrho_{2,i_{00}}\le 0$ for all $i\in a_{00}\left(\bar x,\tilde y\right)$.
Overall, it follows $\tilde i=\mbox{argmax}\left\{c_i\left\vert i\in I_0\left(\bar x\right)\right.\right\}$. According to Lemma \ref{lem:a01}, there have to be $n-s-1$ elements in $\mathcal{E}(\tilde y)$ and we have $\mathcal{E}(\tilde y)\subset {I_0}\left(\bar x\right)\backslash\left\{\tilde i\right\}$ due to the choice of $\tilde i$. This leads to the conclusion that $\tilde y$ could have been constructed as $\bar y$, a contradiction. 

We show that each of the constructed T-stationary points $(\bar x, \bar y)$ of $\mathcal{R}$ is nondegenerate if $\bar x$ has been a nondegenerate M-stationary point of CCOP. 
CC-LICQ at $\bar x$ provides MPOC-LICQ at $(\bar x, \bar y)$ in view of Theorem \ref{thm:cc-mpoc-licq}. Hence, the multipliers for $(\bar x, \bar y)$ defined above are unique. By virtue of NDM2, we have $\bar \mu_q >0$ for all $q \in Q_0(\bar x)$. Moreover, since the components of $c$ are assumed to be positive and pairwise different, we have:
\[
\bar \mu_{2,i}=c_{\bar i}-c_i > 0,i \in \mathcal{\bar E}, \quad \bar \mu_3=c_{\bar i} >0.
\]
Thus, NDT2 is shown.
In case of $\left\|\bar x\right\|_0<s$, NDM3  provides: 
\[
\bar \varrho_{1,i_{00}}=\bar\gamma_{i_{00}} \not =0, i_{00}\in I_0\left(\bar x\right)\backslash\left(\left\{\bar i\right\}\cup\mathcal{\bar E}\right). 
\]
By using Assumption \ref{ass:c} once again, we also have:
\[
\bar \varrho_{2,i_{00}}=c_{i_{00}}-c_{\bar i} < 0, i_{00}\in I_0\left(\bar x\right)\backslash\left(\mathcal{\bar E}\cup \left\{\bar i\right\}\right).
\]
Hence, NDT3 is fulfilled whenever $\left\|\bar x\right\|_0<s$. If the sparsity constrained is active, $a_{00}(\bar x,\bar y)=\emptyset$ and NDT3 trivially holds. 
It is left to show that NDT4 holds. Due to MPOC-LICQ, we immediately obtain the following representation of the tangential space:
\[
\begin{array}{rcl}
        \mathcal{T}^\mathcal{R}_{(\bar x,\bar y)}

&=&\left\{
\xi \in \R^{2n}\,\left\vert\, \begin{array}{l} \begin{pmatrix}

Dh_p(\bar x),0\end{pmatrix} \xi=0, p \in P, \begin{pmatrix}
Dg_q(\bar x),0\end{pmatrix}\xi=0,q \in Q_0(\bar x),\\
\xi_i=0, i \in I_0(\bar x),
\xi_i=0, i \in \{n+1,\ldots,2n\}
\end{array}
\right.\right\}
\end{array}.
\]
Therefore, $D^2 L^\mathcal{R}(\bar x,\bar y)\restriction_{\mathcal{T}^\mathcal{R}_{(\bar x,\bar y)}}$ is nonsingular if and only if $D^2 L(\bar x)\restriction_{\mathcal{T}_{\bar x}}$ is nonsingular.
The latter holds due to NDM4, consequently, NDT4 is also fulfilled. Altogether, the T-stationary point $(\bar x, \bar y)$ is shown to be nondegenerate.
It remains to prove that the T-index of $(\bar x, \bar y)$ equals the M-index of $\bar x$.
The above representation of $\mathcal{T}^\mathcal{R}_{(\bar x,\bar y)}$
 tells us that the number of negative eigenvalues of
$D^2 L^\mathcal{R}(\bar x,\bar y)\restriction_{\mathcal{T}^\mathcal{R}_{(\bar x,\bar y)}}$ is the same as that of $D^2 L(\bar x)\restriction_{\mathcal{T}_{\bar x}}$, i.e. the quadratic indices of $(\bar x, \bar y)$ and $\bar x$ coincide. In order to show that the biactive index $BI$ of $(\bar x, \bar y)$ equals to the sparsity index $SI$ of $\bar x$, we first note:
\[
\left|a_{00}(\bar x, \bar y)\right|=n-\left|a_{10}(\bar x, \bar y)\right|-\left|a_{01}(\bar x, \bar y)\right|. 
\]
By definition, $\left|a_{10}(\bar x, \bar y)\right|=\left\|\bar x\right\|_0$, and by construction, $\left|a_{01}(\bar x, \bar y)\right|=n-s$.
Therefore, we obtain:
\[
BI=\left|a_{00}(\bar x, \bar y)\right|=n-\left\|\bar x\right\|_0-(n-s)=s-\left\|\bar x\right\|_0=SI.
\]
Now, we turn our attention to the additional property NDT5. If $a_{00}\left(\bar x,\bar y\right)\not = \emptyset$, then as above:
\[
\left\|\bar x\right\|_0 = \left|a_{10}(\bar x, \bar y)\right| = n - \left|a_{00}(\bar x, \bar y)\right| - \left|a_{01}(\bar x, \bar y)\right| = n - \left|a_{00}(\bar x, \bar y)\right| - (n -s) = 
s - \left|a_{00}(\bar x, \bar y)\right| < s. 
\]
Hence, the conclusion from NDM3 can be applied and we obtain by recalling $a_{01}(\bar x, \bar y)=\mathcal{\bar E}\cup \left\{\bar i\right\}$:
\[
\bar \sigma_{1,i_{01}}=\bar \gamma_{i_{01}}\ne 0, i_{01}\in a_{01}(\bar x, \bar y).
\]
\qed

The reverse implication of Theorem \ref{thm:mtot} is also valid. From a nondegenerate T-stationary point of $\mathcal{R}$ it can be likewise concluded that the corresponding M-stationary point of CCOP is of the same type.

\begin{theorem}[Stationarity from $\mathcal{R}$ to CCOP]
\label{thm:ttom}
If $(\bar x,\bar y)$ is a T-stationary point of $\mathcal{R}$, then $\bar x$ is an M-stationary point of CCOP. 
If $(\bar x,\bar y)$ is additionally nondegenerate with T-index $m$ and satisfies NDT5, 
then $\bar x$ is also nondegenerate with M-index $m$. 
\end{theorem}
\proof
Since $(\bar x,\bar y)$ is a T-stationary point, there exist
$(\bar \lambda, \bar \mu, \bar \sigma, \bar \varrho)$
with (\ref{eq:tstat-1})--(\ref{eq:tstat-3}). 
The first $n$ rows of (\ref{eq:tstat-1}) read as
\[
   \begin{array}{rcl}
    \nabla f(\bar x)
&=& \displaystyle
\sum\limits_{p \in P}\bar \lambda_p
\nabla h_p(\bar x)
+
    \sum\limits_{q \in Q_0(\bar x)}     \bar \mu_{1,q}
     \nabla g_q(\bar x)+\sum\limits_{i_{01} \in a_{01}\left(\bar x,\bar y\right)} \bar \sigma_{1,i_{01}}
        \cdot e_{i_{01}}
    +\sum\limits_{i_{00} \in a_{00}\left(\bar x,\bar y\right)} \bar \varrho_{1,i_{00}}
        \cdot e_{i_{00}},\\ \\
&=& \displaystyle
\sum\limits_{p \in P}\bar \lambda_p
\nabla h_p(\bar x)
+
    \sum\limits_{q \in Q_0(\bar x)}     \bar \mu_{q}
     \nabla g_q(\bar x)+\sum\limits_{i\in I_0\left(\bar x\right)}
     \bar \gamma_i\cdot e_i,
        \end{array}
\]
where we set:
\[
\bar \mu_{q}=\bar \mu_{1,q}, q \in Q_0(\bar x), \quad \bar \gamma_i=\left\{
\begin{array}{ll}
     \bar \sigma_{1,i},&\mbox{for }i \in a_{01}\left(\bar x,\bar y\right),\\
     \bar \varrho_{1,i},&\mbox{for }i\in a_{00}\left(\bar x,\bar y\right).
\end{array}\right.
\]
It follows that $\bar x$ is an M-stationary point.

Let $(\bar x, \bar y)$ be additionally nondegenerate and fulfill NDT5. First, we note that MPOC-LICQ at $(\bar x, \bar y)$ implies CC-LICQ at $\bar x$ in view of Theorem \ref{thm:cc-mpoc-licq}. NDM2, i.\,e. $\bar \mu_q > 0$ for all $q\in Q_0\left(\bar x\right)$, is also fulfilled as an immediate consequence of NDT2. In order to show NDM3, we assume that the sparsity constraint is not active, i.e. $\left\|\bar x\right\|_0 < s$. Due to Lemma \ref{lem:a01}c), $a_{00}\left(\bar x,\bar y\right)\not =\emptyset$ must hold. NDT5 is applicable and we get $\bar \sigma_{1,i_{01}}\ne 0$ for all $i_{01}\in a_{01}(\bar x, \bar y)$. NDT3 also provides $\bar\varrho_{1,i_{00}}\ne 0$ for all $i_{00}\in a_{00}\left(\bar x,\bar y\right)$. Recalling $I_0(\bar x)= a_{01}\left(\bar x,\bar y\right)\cup a_{00}\left(\bar x,\bar y\right)$ and the definition of $\bar \gamma_i$, $i \in I_0(\bar x)$, NDM3 immediately follows. Finally, we have that NDM4 coincides with NDT4. This follows exactly as in the proof of Theorem \ref{thm:mtot}. From there it can be also seen that the quadratic indices of $\bar x$ and $(\bar x, \bar y)$ coincide. Moreover, for the sparsity index $SI$ of $\bar x$ and for the biactive index $BI$ of $(\bar x, \bar y)$ it is possible to analogously prove $SI=BI$, by using Lemma \ref{lem:a01}c) where needed. Overall, the M-index of $\bar x$ coincides with the T-index of $(\bar x, \bar y)$.
 \qed

%
%

It turns out that, if condition NDT5 is violated, the assertion on the nondegeneracy in Theorem \ref{thm:ttom} does not necessary hold. This becomes clear from the following Example \ref{ex:ndt5}.

\begin{example}[Theorem \ref{thm:ttom} does not hold without NDT5]
\label{ex:ndt5}
We consider the regularized continuous reformulation of CCOP with $n=2$, $s=1$, and $c=(c_1,c_2)^T$ fulfilling $0<c_1<c_2$:
\[
\begin{array}{rl}
\mathcal{R}: \quad \min\limits_{x,y}&(x_1-1)^2+x_2^2+c_1\cdot y_1 + c_2\cdot y_2\\
\mbox{s.t.}&y_1+y_2\ge 1, \quad  
x_i y_i=0,\quad 0\leq y_i\le 1+\varepsilon, \quad i=1,2.
\end{array}
\]
We show that, although $(\bar x,\bar y)=(0,0,0,1) $ is a nondegenerate T-stationary point for $\mathcal{R}$, the point $\bar x=(0,0)$ is degenerate as an M-stationary point for the corresponding CCOP.
We start with the former, where we have $a_{00}(\bar x,\bar y)=\{1\}$, $a_{01}(\bar x,\bar y)=\{2\}$, $a_{10}(\bar x,\bar y)=\emptyset$, $\mathcal{E}(\bar y)=\emptyset$, $\bar y_1+ \bar y_2=n-s=1$.
It holds for the T-stationarity condition:
\[
\begin{array}{rcl}
  \begin{pmatrix}
   -2\\
   0\\
   c_1\\
   c_2
  \end{pmatrix}&=&
\bar \mu_3 \cdot\begin{pmatrix}
0\\
0\\
1\\
1
\end{pmatrix} +\bar \sigma_{1,2}\cdot \begin{pmatrix}
    0\\
    1\\
    0\\
    0
    \end{pmatrix}+\bar \varrho_{1,1}\cdot
    \begin{pmatrix}
   1\\
   0\\
    0\\
    0
    \end{pmatrix}
   + \bar \varrho_{2,1}\cdot
    \begin{pmatrix}
    0\\
    0\\
    1\\
    0    \end{pmatrix}, \end{array}
\]
with the multipliers
\[
\bar \mu_3=c_2>0,\quad \bar \sigma_{1,2}=0,\quad \bar \varrho_{1,1}=-2\ne 0,\quad \bar \varrho_{2,1}=c_1-c_2<0.
\]
Thus, $(0,0,0,1)$ is a T-stationary point fulfilling NDT1-NDT4, but not NDT5.
For $\bar x=(0,0)$ we further have:
\[
\begin{pmatrix}
-2\\0
\end{pmatrix}=
\bar \gamma_1\cdot\begin{pmatrix}
1\\0
\end{pmatrix}
+\bar \gamma_2 \cdot \begin{pmatrix}
0\\1
\end{pmatrix}
\] 
with the unique multipliers $\bar \gamma_1=-2, \bar \gamma_2=0$.
Thus, $\bar x$ is an M-stationary point of CCOP, but it is degenerate due to the violation of NDM3. \qed

\end{example}


Let us discuss the consequences of Theorems \ref{thm:cc-mpoc-licq}, \ref{thm:mtot} and \ref{thm:ttom} for the regularized continuous reformulation $\mathcal{R}$ from the global optimization perspective.
We start by looking at the generic properties of $\mathcal{R}$. 
Here, genericity refers to the fact that
a property of $\mathcal{R}$ holds on an open and dense subset of defining functions $f$, $h$, and $g$ with respect to the strong (or Whitney-) topology, cf. \cite{Jongen:2000}.

\begin{remark}[Genericity for $\mathcal{R}$]
\label{rem:generic}
    We recall from \cite{laemmel:ccop} that CC-LICQ holds generically on the CCOP feasible set. From Theorem \ref{thm:cc-mpoc-licq} we deduce that MPOC-LICQ also holds on the feasible set of a generic regularized continuous reformulation $\mathcal{R}$. Moreover, it is known that all M-stationary points of CCOP are generically nondegenerate, see \cite{laemmel:ccop}. In view of  Theorem \ref{thm:ttom}, any T-stationary point of $\mathcal{R}$ is induced by an M-stationary point of CCOP.
    Hence, Theorem \ref{thm:mtot} provides that all T-stationary points of $\mathcal{R}$ are also generically nondegenerate and additionally fulfill NDT5. This suggests that the assumption of NDT5 is not restrictive. \qed
\end{remark}


We turn our attention to the Morse theory for the regularized continuous reformulation $\mathcal{R}$.

\begin{remark}[Morse theory for $\mathcal{R}$]
\label{rem:morse}
Since we know by Remark \ref{rem:generic} that all T-stationary points of $\mathcal{R}$ are generically nondegenerate, the Morse theory developed for MPOC in \cite{laemmel:mpoc} can be successively applied for its subclass $\mathcal{R}$. Morse theory typically provides deformation and cell-attachment results for the lower level sets of an underlying optimization problem. Deformation means in our context that outside the set of T-stationary points the topology of lower level sets of $\mathcal{R}$ remains unchanged if the level varies. Cell-attachment describes the topological changes of the lower
level sets if passing a nondegenerate T-stationary point. 
Namely, a cell of dimension equal to the T-index needs to be attached to the lower level set. By 
doing so, we get the lower level set corresponding to the T-stationary point up to a homotopy. Based on deformation and cell-attachment results, it is possible to relate the numbers of minimizers and saddle points, thus, to adequately capture the global structure of $\mathcal{R}$. E.g., the so-called mountain pass result from global optimization says that there exist at least $(k-1)$ T-stationary points of $\mathcal{R}$ with T-index equal to one, where $k$ denotes the number of local minimizers of $\mathcal{R}$. These saddle points can be of two types: either (I) $BI=0$ and $QI=1$ or (II) $BI=1$ and $QI=0$. Whereas the saddle points of type (I) appear quite similarly in the nonlinear programming, type (II) is characteristic for the regularization $\mathcal{R}$. In particular, at the saddle points of type (II) the biactive index set does not vanish. This phenomenon cannot thus be neglected if studying T-stationary points of $\mathcal{R}$. \qed
\end{remark}

For the global structure of $\mathcal{R}$, it is valuable to make a comparison to CCOP concerning the numbers of their nondegenerate M- and T-stationary points with the same M- and T-index, respectively. Let us start by considering minimizers. 


\begin{corollary}[Minimizers of CCOP and of $\mathcal{R}$]
\label{cor:min}
It holds:
\begin{itemize}
    \item[a)] If $\bar x$ is a nondegenerate minimizer of CCOP, then there exists unique $\bar y$, such that $(\bar x, \bar y)$ is a minimizer of $\mathcal{R}$, moreover, it is nondegenerate. 
    \item[b)] If $(\bar x, \bar y)$ is a nondegenrate minimizer of $\mathcal{R}$, then $\bar x$ is a nondegenerate minimizer of CCOP.
\end{itemize}

\end{corollary}

\proof 
Lemma 3 from \cite{laemmel:ccop} says that a nondegenerate minimizer $\bar x$ of CCOP has zero M-index and, in particular, $\left\|\bar x\right\|_0=s$. In view of Theorem \ref{thm:mtot},
we immediately obtain the assertion. For the reverse implication, we apply Corollary 1 from \cite{laemmel:mpoc} saying that the T-index of a nondegenerate minimizer $(\bar x,\bar y)$ of an MPOC must vanish. In particular, it holds $a_{00}(\bar x, \bar y)= \emptyset$. Then, NDT5 is trivially satisfied, and by Theorem \ref{thm:ttom} we are done. 
\qed

From Corollary \ref{cor:min} we conclude that the numbers of minimizers of CCOP and $\mathcal{R}$ coincide, at least in a generic situation where they are nondegenerate. However, the numbers of saddle points with nonvanishing M- and T-index, respectively, differ. 
 

\begin{remark}[Saddle points of CCOP and  of $\mathcal{R}$]
\label{rem:saddle}
    If at an nondegenerate M-stationary point $\bar x$ the sparsity constraint is not active, i.e. $\left\|\bar x\right\|_0 < s$, it induces multiple T-stationary points $(\bar x, \bar y)$ of $\mathcal{R}$. Due to Theorem \ref{thm:mtot}, their precise number is $\binom{n-\left\|\bar x\right\|_0-1}{n-s-1}$. Moreover, the quadratic indices of $\bar x$ and $(\bar x, \bar y)$ coincide. The sparsity index $SI=s-\left\|\bar x\right\|_0$ of $\bar x$ and the biactive index $BI=\left|a_{00}(\bar x, \bar y)\right|$ of $(\bar x, \bar y)$ are equal. The appearance of exponentially many T-stationary points of the same type can be explained in terms of the cell-attachment for CCOP. In \cite{laemmel:ccop}, it is shown that for describing topological changes of the CCOP lower level sets the attachment of multiple cells of dimension equal to the M-index is needed. The number of cells to be attached amounts to $\binom{n-\left\|\bar x\right\|_0-1}{s-\left\|\bar x\right\|_0}$. However, it holds:
    \[
        \binom{n-\left\|\bar x\right\|_0-1}{n-s-1} = 
        \binom{n-\left\|\bar x\right\|_0-1}{n-\left\|\bar x\right\|_0-1-(n-s-1)} = \binom{n-\left\|\bar x\right\|_0-1}{s-\left\|\bar x\right\|_0}.
    \]
   We see that the number of induced T-stationary points of $\mathcal{R}$ corresponds to the number of cells to be attached if passing the corresponding M-stationary level in CCOP. This observation suggests that the involved structure of saddle points in CCOP translates into the increasing number of saddle points in $\mathcal{R}$.  \qed  
\end{remark}


\section{Comparison to the literature}

At the end of this section we compare our results on the regularized continuous reformulation $\mathcal{R}$ with those achieved in the literature for the original continuous reformulation (\ref{eq:relax}). Note that $\mathcal{R}$ becomes (\ref{eq:relax}) if we set $c=0$ and $\varepsilon=0$. In \cite{burdakov:2016}, it has been shown how minimizers of CCOP and of (\ref{eq:relax}) are related.

\begin{proposition}[Minimizers of CCOP and of (\ref{eq:relax}), see \cite{burdakov:2016}]
    \label{cor:min-old}
It holds:
\begin{itemize}
    \item[a)] If $\bar x$ is a minimizer of CCOP, then there exists $\bar y$, such that $(\bar x, \bar y)$ is a minimizer of (\ref{eq:relax}). If additionally, $\left\|\bar x\right\|_0=s$, then $\bar y$ is unique.
    \item[b)] If $(\bar x, \bar y)$ is a minimizer of (\ref{eq:relax}) with $\left\|\bar x\right\|_0=s$, then $\bar x$ is a minimizer of CCOP.
\end{itemize}
\end{proposition}

\noindent
From here we see that there is a one-to-one correspondence between the minimizers of CCOP and of (\ref{eq:relax}), whenever the sparsity constraint is active. A similar result follows for the regularized continuous reformulation $\mathcal{R}$ from Corollary \ref{cor:min}, since at a nondegenerate minimizer of CCOP the cardinality constraint is active. Concerning the general relation between the M-stationary points of CCOP and the T-stationary points of (\ref{eq:relax}), the following has been shown in \cite{laemmel:mpoc}.

\begin{proposition}[Stationarity relations between CCOP and (\ref{eq:relax}), see \cite{laemmel:mpoc}]
    A feasible point $(\bar x, \bar y)$  of (\ref{eq:relax}) is T-stationary if and only if the point $\bar x$ is M-stationary for CCOP.
\end{proposition}

\noindent
Although this appealing relation holds for (\ref{eq:relax}), cf. a similar result for $\mathcal{R}$ in Theorems \ref{thm:mtot} and \ref{thm:ttom}, the original continuous reformulation is intrinsically degenerate.

\begin{proposition}[Degeneracy of (\ref{eq:relax}), see \cite{laemmel:mpoc}]
\label{prop:deg}
All T\Hyphdash stationary points of  (\ref{eq:relax}) are degenerate.
\end{proposition}

\noindent
The degeneracy of T-stationary points of (\ref{eq:relax}) prevails even when we start by a CCOP with all its M-stationary points being nondegenerate. We emphasize that this hampers the study of (\ref{eq:relax}) from the global optimization perspective. In particular, it is not possible to develop the Morse theory for (\ref{eq:relax}) and, finally, to conclude that the price to pay for reformulating the cardinality constraint is the appearance of exponentially many additional saddle points, cf. Remark \ref{rem:saddle}. Actually, the saddle points of (\ref{eq:relax}) may well constitute a continuum, rather than to be isolated as it is generically the case for $\mathcal{R}$.

\begin{example}[Continuum of T-stationary points, see \cite{laemmel:mpoc}]
  We consider the following CCOP:
  \[
\min_{x \in \R^2}\,\, \left(x_1-1\right)^2+\left(x_2-1\right)^2\quad \mbox{s.\,t.} \quad \NNorm{x}{0} \le 1.
\]
Its M-stationary points are $(1,0)$, $(0,1)$, and $(0,0)$, all of them being nondegenerate.
The corresponding continuous reformulation (\ref{eq:relax}) reads as
\[
   \min_{x,y \in \R^2} \,\,\left(x_1-1\right)^2+\left(x_2-1\right)^2 \quad \mbox{s.\,t.} \quad 
   y_1+y_2 \geq 1,  \quad 
      y_1,y_2 \in [0,1],  \quad 
      x_1y_1 =0, x_2y_2=0.
\]
We get as its T-stationary points:
\[
(1,0,0,1), (0,1,1,0), \mbox{ and } (0,0,y_1,y_2) \mbox{ with } y_1, y_2\in [0,1], y_1+y_2 \ge 1.
\]
By Theorem \ref{prop:deg} all points are degenerate. Moreover, none of the saddle points is isolated.
\qed
\end{example}

\section*{Conclusions}

We emphasize that the introduction of auxiliary $y$-variables shifts the
complexity of dealing with the cardinality constraint in CCOP into the appearance of multiple saddle points for its continuous reformulation. For this conclusion to make, we appropriately regularized the original continuous reformulation (\ref{eq:relax}). The regularized continuous reformulation $\mathcal{R}$ turns out to have the same favorable properties in regard to the minimizers of CCOP. In addition, we could successively apply the Morse theory developped for MPOC in \cite{laemmel:mpoc} to $\mathcal{R}$. This is possible due to the generic nondegeneracy of its T-stationary points, whereas the T-stationary points of (\ref{eq:relax}) are intrinsically degenerate. The global structure of the (regularized) continuous reformulation is therefore fully understood. Next step would be to study the global structure of the so-called Scholtes-type regularization of the continuous reformulation along the same lines. This is the topic of current research.   

\bibliographystyle{plain}
\bibliography{lit.bib}
\end{document}